\def\obs#1{{\bf (*** #1 ***)} }

 \def\obs#1{}     
 
\documentclass[twoside,letterpaper,draft,11pt]{amsart}

\usepackage{amsmath,amsthm,latexsym,xspace,amscd,amssymb,xypic,amscd,amssymb}               
\usepackage{color}
\usepackage{fullpage}

\newtheorem{teo}{Theorem}[section]
\newtheorem{defi}[teo]{Definition}
\newtheorem{lema}[teo]{Lemma}

\newtheorem{prop}[teo]{Proposition}

\newtheorem{exe}[teo]{Example}

\newcommand{\X}{{\mathbb X}}
\newcommand{\Y}{{\mathbb Y}}

\newcommand{\Z}{{\mathbb Z}}

\newcommand{\ex}{{\exists}}
\usepackage[usenames,dvipsnames]{xcolor}
\newcommand{\m}{{}^{-1}}
\newcommand{\mt}{\mapsto}

\def\ndv{\ {\mid \kern -0.7 em {\scriptstyle \not}} \ \ }
\def\nd{\ {\mid \kern -0.4 em {\scriptstyle \not}} \ \ }

\newcommand{\N}{{\mathbb N}}

\numberwithin{equation}{section}
\title[Borel Globalization]{Borel globalizations of partial actions of  Polish groups}

\begin{document}

\author[H.\ Pinedo ]{H. Pinedo }
\address{Escuela de Matem\'aticas, Universidad Industrial de Santander, Cra. 27 Calle 9 UIS Edificio 45\\  Bucaramanga, Colombia}\email{hpinedot@uis.edu.co }
\author[C.\ Uzcategui ]{C. Uzcategui }
\address{Escuela de Matem\'aticas, Universidad Industrial de Santander, Cra. 27 Calle 9 UIS Edificio 45\\  Bucaramanga, Colombia}\email{cuzcatea@uis.edu.co}
\date{\today}

\begin{abstract}  We show that the enveloping space $\X_G$ of a partial action of a Polish group $G$ on a Polish space $\X$ is a standard Borel space, that is to say,  there is a topology $\tau$ on $\X_G$ such that $(\X_G, \tau)$ is Polish and  the quotient Borel structure on $\X_G$ is equal to $Borel(\X_G,\tau)$. To prove this result we show a generalization of a theorem of Burgess about Borel selectors for the orbit equivalence relation induced by a group action and also show that some properties of the Vaught's transform are valid for partial actions of groups. 

\end{abstract}

\keywords{Partial action, enveloping space, Polish space, standard Borel space. }
\subjclass[msc2010]{54H15, 54E50, 54E35}

\maketitle
\section{Introduction}

Given an action $a:G\times \Y\rightarrow \Y$ of a group $G$ over a set $\Y$ and an invariant  subset $\mathbb{X}$ of $\Y$ (i.e. $a(g,x)\in \mathbb{X},$ for all $x\in \mathbb{X}$ and $g\in G$),  the restriction of $a$ to $G\times\mathbb{X}$ is an action of $G$ over $\mathbb{X}$. However, if $\mathbb{X}$ is not invariant,  we get what is called a {\em partial action} on $\mathbb{X}$: a collection of partial maps $\{m_g\}_{g\in G}$ on  $\mathbb{X}$ for which  $m_1={\rm id}_\mathbb{X}$ and $m_g\circ m_h\subseteq m_{gh},$ for all $g,h\in G$.
A partial action of a group   is a weakening of the classical notion of a group action and was introduced by  R. Exel in \cite{E1} motivated by  problems arising from $C^*$-algebras (see \cite{E-1},  \cite{E-2},  \cite{E-3},  \cite{E0},  \cite{E2}).  Since then, partial group actions  have appeared in many different context, such as the  theory of operator algebras and in algebra  \cite{D2}, also in the theory of $\mathbb{R}$-trees, tilings and model theory  \cite{KL}. More recently partial actions have been an efficient tool to study $C^*$-\'algebras associated with integral domains   \cite{BoE1} and dynamical systems of type  
 $(m,n)$ \cite{AraEKa}.  A deep application of partial actions   was obtained in   \cite{AraE}  in which the authors provide a counter-example  for a  conjecture of M. R{\o}rdam and  A. Sierakovski related with the  Banach-Tarski paradox.

 In the topological context, partial actions on topological spaces, consists of a family of homeomorphism between open subsets of the space (see  Definition \ref{topar}).
Extension problems in topology are related to  the question of  whether a given collection of partial maps on a space can be realized as the set of traces of a corresponding
collection of total maps on some superspace. In the context of partial actions this question was first  studied in \cite{AB} and \cite{KL}, where the authors showed that every partial action  of a topological group $G$ on a topological space $\X$ can be obtained as a restriction of a continuous global action of $G$ on some superspace $\X_G$ of $\X,$ such $\X_G$ is called the globalization or enveloping space of $\X$ (see  Definition \ref{globa}). Similar approaches were also presented in \cite{Megre2} and \cite{St2}.

Actions of Polish groups have received a lot of attention  in recent  years, because of  its connection with many areas of mathematics (see \cite{KEB,GA} and the references therein).
We recall that   a Polish space  is a topological space which is separable and completely metrizable, and a Polish group is a topological group whose topology is Polish. The problem  of whether a partial action of a Polish group on a Polish space admits a Polish globalization was studied in \cite{PU}.  In general, the spaces $\X_G$ and $\X$ do not share the same topological properties. For example, $\X_G$ could be  not Hausdorff even for $\X$ and $G$ metric spaces  (see for instance  \cite[Example 1.3, Proposition 1.2]{AB} and \cite[Proposition 2.1]{EG}). It was shown in \cite{PU} that under appropriate hypothesis $\X_G$ is a Polish space (see Theorem \ref{xgpol} below).   

On the other hand, the classical theory of Polish group actions has  been  extended to Borel actions of Polish groups, where the action is required to be  Borel measurable rather than   continuous, an approach that goes back to the works of Mackey \cite{Mackey1957} and Effros \cite{EF} and many others (see for instance \cite{KEB, GA} for a more recent treatment).  In this article we present a first step for  extending the notion of a topological globalization to the category of standard Borel spaces.  One of our  main results is that, when $\X$ and $G$ are Polish spaces, the enveloping space $\X_G$ is a standard  Borel space,  that is there is a topology $\tau$ on $\X_G$ extending the topology of $\X_G$ such that $(\X_G, \tau)$ is Polish and  the quotient Borel structure of $\X_G$ is equal to  $Borel(\X_G,\tau)$. Moreover, the Borel structure of $\X$ is not changed, that is,  $Borel(\X)=Borel(\X,\tau_\X)$ and the global action of $G$ over $\X_G$ is Borel measurable (see Theorem \ref{gloto}). In other words, $(\X_G,\tau)$ is a Borel $G$-space (as in \cite{KEB}) extending the original partial action on $\X$ and thus it can be called a Borel globalization of $\X$.

For the proof of our results we need to extend  some well known theorems about Polish group actions to the context of partial actions. In particular, we present a generalization of a theorem of Burgess \cite{Bu79}  about Borel selectors for the orbit equivalence relation  $E^p_G$ on $\X$ induced by the partial action of a Polish group $G$.  In addition, we show some properties of   the natural generalization of the  Vaught's transforms \cite{vaught74} but now for partial actions. Vaught's transforms  are an important tool in invariant descriptive  set theory \cite{GA,KE}, and they have also  been studied in the context of Polish groupoid  \cite{LU}.   Since $\X_G$ is a standard Borel space, then it is natural to look at the orbit equivalence relation  $E_G$ associated to the enveloping action of $G$ over $\X_G$. We show that the orbit equivalence relation $E^p_G$ given by the partial action is Borel bireducible to  $E_G$.

\section{Preliminaries}
\label{not}

Let $G$ be a  group with identity element $1$, $\X$ a set and  $m\colon  G\times \X\to \X,\,\,(g,x)\mt m(g,x)=g\cdot x\in \X$  a partially defined function.
We write $\ex g\cdot x$ to mean that $(g,x)$ is in the domain of $m.$ Then $m$  is called a (set theoretic) {\it partial action} of $G$ on $\X,$ if for all $g,h\in G$ and $x\in \X$ we have:
\smallskip

\noindent (PA1) If $\ex g\cdot x$, then  $\ex g\m\cdot (g\cdot x)$ and $g\m\cdot (g\cdot x)=x$,
\smallskip

\noindent (PA2)  If $\ex g \cdot (h\cdot x)$, then  $\ex (g h)\cdot x$ and  $g \cdot (h\cdot x)=(g h)\cdot x$,
\smallskip

\noindent (PA3) $\ex 1\cdot x,$ and $1\cdot x=x$.
 \smallskip

We fix some notations that shall be useful throughout the work. Let 
$$
G* \X=\{(g,x)\in G\times \X\mid \ex g\cdot x \}.
$$
Moreover, for $h\in G$ and $y\in \X$ write
$$
\X_{h\m}=\{ x\in \X\mid \ex h\cdot x\},
$$
 and $m_h\colon \X_{h\m}\ni x\mt h\cdot x\in\X_h.$   Finally set $G^y=\{g\in G\mid \ex g\cdot y\},$   $G_y=\{g\in G^y\mid g\cdot y=y\},$ the {\it stabilizer} of $y$ and  $G^y\cdot y=\{g\cdot y\mid g\in G^y\},$ the {\it orbit} of $y$.


By    \cite[Lemma 1.2]{QR} a partial action $m\colon G*\X\to\X$  can be equivalently  formulated in terms  a family of bijections  in the following sense.
 A partial action $m$ of $G$ on $\X$ is a family $m=\{m_g\colon \X_{g\m}\to\X_g\}_{g\in G},$ where $\X_g\subseteq \X,$  $m_g$ is bijective, for  all $g\in G,$  and such that:
\begin{itemize}
\item[(i)]$\X_1=\X$ and $m_1=\rm{id}_\X;$
\item[(ii)]  $m_g( \X_{g\m}\cap \X_h)=\X_g\cap \X_{gh};$
\item[(iii)] $m_gm_h\colon \X_{h\m}\cap  \X_{ h\m g\m}\to \X_g\cap \X_{gh},$ and $m_gm_h=m_{gh}$ in $ \X_{h\m}\cap  \X_{g\m h\m};$
\end{itemize}
for all $g,h\in G.$

\begin{exe}
\label{indu} {\bf  Induced partial action:}
Let $u \colon G\times \Y\to \Y$ be an action of $G$ on $\Y$ and $\X\subseteq \Y.$ For $g\in G,$ set $\X_g=\X\cap u_g(\X)$ and let $m_g=u_g\restriction \X_{g\m}.$  Then $m\colon G* \X\ni (g,x)\mt m_g(x)\in \X $ is a partial action of $G$ on $\X.$  In this case we say that $m$ is induced by $u.$
\end{exe}

From now on $G$ will be a topological group and $\X$ a topological space. We consider the set  $G\times \X$ with the product topology and    $G* \X\subseteq G\times \X$  with the subspace topology. 

\begin{defi}\label{topar} A topological partial action of the group $G$ on the topological space $\X$ is a partial action $m=\{m_g\colon \X_{g\m}\to\X_g\}_{g\in G}$ on the underlying
set $\X$ such that each $\X_g$ is open in $\X$, and each $m_g$ is a homeomorphism.
\end{defi}

We recall the definition of the enveloping action in the topological sense. Let $m$ be a topological partial action of $G$ on $\X.$ Define an
equivalence relation on $ G\times \X$ as follows:
\begin{equation}
\label{eqgl}
(g,x)R  (h,y) \Longleftrightarrow x\in \X_{g\m h}\,\,\,\, \text{and}\, \,\,\, m_{h\m g}(x)=y,
\end{equation}
and denote by  $[g,x]$   the equivalence class of the pair $(g,x).$ Consider the set $\X_G=(G\times \X)/R$  endowed with the quotient topology.  Then by  \cite[Theorem 1.1]{AB}  the action 
\begin{equation}
\label{action}
\mu \colon G\times \X_G\ni (g,[h,x])\to [gh,x]\in \X_G,
\end{equation}
is continuous and the map
\begin{equation}
\label{iota}
\iota \colon \X\ni x\mt [1,x]\in \X_G
\end{equation} is a continuous injection.

The following is the basic result about  topological partial actions.

\begin{teo}
\label{Abadie} \cite[Theorem 1.1]{AB} and  \cite[Propositions 3.11 and  3.12]{KL}. Let $m$ be a continuous partial action of $G$ on $\mathbb{X}$. Then
\begin{enumerate}
\item[(i)]   $\iota : \mathbb{X}\rightarrow  \iota(\mathbb{X})$ is a homeomorphism. In addition, if $G*\X$ is open in $G\times \X$, then $\iota(\X)$ is open in $\X_G$.
\item[(ii)] $m$ is equivalent to the partial action induced on $\iota(\mathbb{X})$ by the global action $\mu,$ that is  $\mu_g(\iota(x))=\iota(m_g(x)),$ for all $g\in G$ and $x\in \X_{g\m}.$
\end{enumerate}
\end{teo}

\begin{defi}\label{globa} Let $m$ be a partial action of $G$ on $\X.$ The action 
provided by \eqref{action} is called the enveloping action of $m$ and  the space $\X_G$ is the
enveloping space  or a globalization of $\X.$ \end{defi}

An important concept in the study of group actions is the associated orbit equivalence relation. Next we introduce the natural generalization of this concept to the context of partial actions. 

\begin{defi}  Let $m$ be a partial action on $\X$. The orbit equivalence relation $E_G^p$ on $\X$ is defined by
$$x E_G^p y \Longleftrightarrow \ex\, g\cdot x\,\,\,\text{and}\,\,\, g\cdot x=y ,$$ for some $g\in G.$ If $m$ is global, we simply write $E_G$ instead of $E_G^p.$
\end{defi}
The equivalence class of  $x\in\X$ is $[x]=G^x\cdot x$. The set of equivalence classes  $\X/E_G^p$  is endowed with the quotient topology. We have the following.

\begin{lema} 
\label{gxg} 
Given a   topological partial action $m=\{m_g\colon \X_{g\m}\to\X_g\}_{g\in G}$  of $G$ on  $\X$. The following assertions hold:
\begin{enumerate}
\item[(i)] For any $g\in G$ and $x\in \X_g$ we have $G^xg=G^{g\m\cdot x}.$ In particular, $G^x$ and $G^y$ are homeomorphic, if $xE_G ^p y$.
\item[(ii)]  The quotient map $\pi \colon \X\ni x\mapsto [x] \in \X/E_G^p$ is continuous open and if  $\X$ is second countable, then so is $\X/E_G^p.$ 
\item[(iii)]  Let $h,g\in G$ and $x\in \X$. Suppose $g\in G^x$  and $g^{-1}h\in G_x$. Then $h\in G^x$.
\end{enumerate}
\end{lema}
\proof (i)   See \cite[Lemma 5.2]{PU}.
 (ii) That $q$ is continuous and open is  Lemma 3.2 of  \cite{PU}. Let $\{U_n\}_{n\in\N}$ be a countable basis of $\X.$ Since $\pi$ is open,   the family $\{\pi[U_n]\}_{n\in\N}$ consists of open subsets of $ \X/E^p_G$ and  the continuity of $\pi$ implies that $\{\pi[U_n]\}_{n\in\N}$ is also a basis.

 (iii) Since $x=(g^{-1}h)\cdot x,$ and $g\in G^x,$  then $g\cdot x=g\cdot( (g^{-1}h))\cdot x=h\cdot x,$ thanks to  (PA2),  and the result follows.
 \endproof
 
The enveloping space $\X_G$ can also be presented as a quotient with respect to an orbit equivalence relation of a partial action. This fact will be crucial for the rest of the paper. 
 
Let  $m$ be a partial action of $G$ on $\X,$ and denote $(G\times \X)_{g}=G\times \X_{g}$,  for any $g\in G.$ Then the family
$
\widehat{m}=\{\widehat{m}_g\colon (G\times \X)_{g\m}\to (G\times \X)_{g}\}_{g\in G},\,\, \widehat{m}_g(h,x)=(hg\m, g\cdot x),
$
is a partial action of $G$ on $G\times \X$. Denote by $\widehat{E}^p_G$ the orbit equivalence relation induced by $\widehat{m}.$  The importance of this partial action $\widehat m$ is the following.

\begin{teo}\cite[Theorem 3.3]{PU}
\label{mtomg}
Let  $m$ be  a topological partial action of $G$ on $\X$, and $R$ be the equivalence relation given by \eqref{eqgl}. Then
$(g,x)\widehat{E}^p_G(h,y)$ iff $ (g, x)R(h,y)$, for all $g,h\in G$ and $x,y\in \X$. Hence, the orbit space $(G\times \X)/\widehat{E}^p_G$ coincides with  $\X_G.$
\end{teo}

A {\em Polish space} is a separable topological space which is metrizable by a complete metric.  We recall that a subspace of Polish space $\X$ is Polish iff it is a $G_\delta$ subset of $\X$ \cite[Theorem 3.11]{KE}.

The following result,  mentioned in the introduction,  indicates under which conditions we can obtain a Polish globalization.

\begin{teo}\label{xgpol}\cite[Theorem 4.7]{PU}
 Let $m$ be a continuous partial action of a separable metrizable group $G$ on a separable metrizable space  $\mathbb{X}$ such that $G*\mathbb{X}$ is open  and $\widehat{E}^p_G$ is closed. If  $G*\mathbb{X}$ is clopen or $\mathbb{X}$ is locally compact, then $\mathbb{X}_G$ is metrizable. If in addition,  $G$ and $\mathbb{X}$ are Polish, then  $\mathbb{X}_G$ is Polish.
\end{teo}

{\sf For the rest of this paper,  $G$ will be a Polish group, $\X$ a Polish space and all partial actions are assumed to be continuous with $G*\X$ is $G_\delta$}. 
\smallskip

 Most of the results in \cite{AB,KL} are stated for partial actions where $G*\X$ is open. However, in \cite[Remark 1.1]{AB} provides  an example of a topological partial action in which $G*\X$ is closed in $G\times \X$.  Abadie's example is a particular case of the following.

\begin{exe} Let $G$ be a Polish group, $H$ a Polish subgroup of $G$ and $a:H\times \X\to \X$ an action of $H$ in a Polish space $\X$. Consider $m=\{m_g\colon \X_{g\m}\times \X_g\}_{g\in G},$ where  $m_g=a_g$ if $g\in H$ and $\X_g=\emptyset,$ otherwise. Then $m$  is a topological partial action of $G$ on $\X$ for which $G*\X=H\times \X$ is $G_\delta.$
\end{exe}

As we said in the introduction, we are interested in the Borel structure of the enveloping space, so we recall that the  the {\it quotient Borel structure} on $\X/E^p_G$   consists of all subsets $B$ of $\X/E^p_G$ such that $\pi\m(B)$ is Borel in $\X.$

Now  we introduce some  terminology about equivalence relations (see \cite{KEB,GA}). Let $E$ and $F$ be equivalence relations on some topological spaces $\X$ and $\Y$, respectively. 
A  {\em reduction} of $E$ into $F$  is a map  $f:\X\rightarrow \Y$ such that $xEy$ iff $f(x)F f(y)$. We say that $E$ and $F$ are {\em Borel bireducible}, if there are reductions of $E$ into $F$ and viceversa which are Borel measurable.   A {\em Borel selector} for  $E$  is a Borel function $S:\X\rightarrow \X$ such that $x Ey$ iff $S(x)=S(y)$ and $S(x)Ex,$ for all $x\in\X$.  We say that $E$ is {\em smooth} if there is a Borel reduction  of $E$ into the identity relation of some Polish space $\Y$, i.e, if there is a Borel map $f:\X\rightarrow \Y$ such that $xEy$ iff $f(x)=f(y)$. If there is a Borel selector for $E$, then $E$ is smooth, however, the reciprocal does not hold in general.  

We will need the following result about $G_\delta$ equivalence relations (i.e. $E$ is $G_\delta$ as a  subset of $\X\times\X$) which is part of  the   Glimm-Effros dichotomy, a fundamental result about Borel equivalence relations  \cite{HKL}. 

\begin{teo}\cite[Theorem 1.1]{HKL}
\label{gdelta}
Any  $G_\delta$ equivalence relation over a Polish space is smooth.
\end{teo}

The $\sigma$-algebra  of Borel subsets of a topological space $\X$ is denoted by $Borel(\X)$. A space $\X$ is a {\em standard Borel space} \cite{KE} if there is a Polish space $\Y$ such that $\X$ and $\Y$ are Borel isomorphic, equivalently, if there is a Polish topology $\tau$ on $\X$ such that $Borel(\X)=Borel(\X,\tau)$.  

\section{The Vaught's transforms and Borel selectors for partial actions}
 
Let  $m$ be a topological partial action of  $G$ on $\X.$ Given $x\in \X$ and $V\subseteq G$ denote $V^x=V\cap G^x.$ Then for  $A\subseteq \X$ and a  nonempty open set $V\subseteq G,$ the Vaught transforms are
$$
A^{\triangle V}=\{ x\in\X\mid \{ g\in V^x \colon  g\cdot x\in A \}\; \mbox{is not meager in $V^x$}\} 
$$
and 
\[
A^{* V}=\{x\in\X\mid \{ g\in V^x\colon  g\cdot x\in A \}\,\mbox{is comeager in $V^x$}\}.
\]
As usual, using the quantifier  $\forall^* g\in O$ which is read as saying  for comeager many $g$ in $O$  and   $\exists^* g\in O $ as for nonmeager many $g$ in $O$, the Vaught's transform are stated as follows:
$$
A^{\triangle V}=\{ x\in\X\mid \exists^* g\in V^x\,\, g\cdot x\in A \}\,\,\,\text{and}\,\,\,A^{* V}=\{x\in\X\mid \forall^* g\in V^x\,\, g\cdot x\in A \}.
$$

Now we show that some properties of the Vaught transforms known to hold for total actions (see \cite[16B, p. 95]{KE}) also hold for partial actions. 
The proof of the following proposition is completely  analogous to the case of a global action.

\begin{prop}
\label{VT} Let  $m$ be  a topological partial action of  $G$ on $\X$. Then  for  $A\subseteq \X$,  a  nonempty open set $V\subseteq G$ and  a family  $\{A_n\}_{n\in\N}$ of subsets of $\X$ the following  assertions hold:
\begin{enumerate}
\item[(i)] $ \X\setminus A^{\triangle V}=(\X\setminus A)^{* V}\,\,\, \text{and} \,\,\,\,\X\setminus A^{* V}=(\X\setminus A)^{\triangle V}$,
\item[(ii)] If   $A=\bigcup_{n} A_n$, then  $ A^{\triangle V}=\bigcup_{n} A_n^{\triangle V}$,
\item[(iii)]  If $A=\bigcap_{n} A_n$, then $A^{*V}=\bigcap_{n} A_n^{*V}$,

\item[(iv)] Let $\{U_n\}_n$ be a basis of open sets of $G$. If $A$  is analytic, then $A^{\triangle V}=\bigcup \{A^{*U_n}:\; U_n\subseteq V\}$.
\end{enumerate}
\end{prop}

\begin{teo}
\label{VT3}
Let  $m$ be a  partial action of $G$ on $\X.$ Let  $V \subseteq G$ a nonempty open set. 
If $A\subseteq \X$ is Borel, then  $A^{\triangle V}$ and $A^{*V}$ are Borel. 
\end{teo}
\proof
We first show that if $A\subseteq \X$ is  open, then $A^{\triangle V}=\bigcup_{g\in V} \{x\in \X_{g^{-1}} :\; g\cdot x\in A \}$ and $A^{\triangle V}$ is open. 
Let $A\subseteq \X$ be  open and  $x\in A^{\triangle V}$. Since $\{g\in V^x \mid g\cdot x\in A \} $ is not meager,  there exists  $g\in  V^x$ such that $g\cdot x\in A.$ In particular $x\in \X_{g^{-1}}$. Conversely, suppose $x\in \X_{g^{-1}}$ and $ g\cdot x\in A$ for some $g\in V$.  
Notice that  $\{h\in V^x:\; h\cdot x\in A\}$ is open in the Polish space  $V^x$, thus  it is either empty or not meager in $V^x$.  
Since  $g\in  \{h\in V^x: h\cdot x\in A\}$, then this  set is not meager in $V^x$. Hence    $x\in A^{\triangle V}$.
The rest follows from Proposition \ref{VT}.\endproof

Now  we present  generalizations to the context of partial actions of some well known facts about selectors for orbit equivalence relations induced by Polish group actions. 

\begin{defi}
Let $\X$ be a standard Borel space and E an equivalence relation on $\X$. We
call E idealistic if there is an map $C \mapsto I_C$  associating to each
E-equivalence class C a $\sigma$-ideal $I_C$ of subsets of C such that:
\begin{itemize}
\item[(i)] $C\notin I_C;$

\item[(ii)]  for each Borel set $A\subseteq X\times X$ the set $A_I$ defined by
\begin{equation*}
\label{bor}
x \in A_I \Leftrightarrow \{y \in [x] : (x, y) \in A\} \in I_{[x]}\end{equation*}
is Borel.
\end{itemize}
\end{defi}

The orbit equivalence relation induced by  a Polish group action is idealistic (see \cite{KE1} and also \cite[Theorem 5.4.10]{GA}). The analogous result holds for partial actions. 

\begin{teo}
\label{idealisitic}
Let $m$ be a partial action of $G$ on $\X$. Then $E^p_G$
 is idealistic.
\end{teo}

\proof
Let $x\in\X$ and set $C=[x]=G^x\cdot x.$ Consider the $\sigma$-ideal  $I_C$   of $\mathcal{P}(C)$ defined by
$$
S\in I_C\Longleftrightarrow \{g \in G^x : g\cdot x \in S\}\,\,\text{is meager in}\,\, G^x.
$$ 
Observe that (i) of Lemma \ref{gxg} guarantees that this definition does not depend on the representative of $[x]$.  Since the set $\{g \in G^x : g\cdot x \in [x]\}=G^x$ is a non empty open set, then $C\notin I_C.$ Now we check  that for any  Borel set  $A\subseteq  \X\times \X$, the set $A_I$ defined in \eqref{bor} is Borel. Indeed
\begin{align*}
x \in A_I& \Longleftrightarrow \{y \in [x] : (x, y) \in A\} \in I_{[x]}\\
&\Longleftrightarrow  \{g \in G^x : (x, g \cdot x) \in A\} \,\,\text{is meager in}\,\, G^x.
\end{align*}
Consider the topological  partial action of $\beta$ of $G$ over $\X\times\X$ given by the family of partial homeomorphisms $\{ \X\times \X_{g\m}\ni  (x,y)\to (x, g\cdot y)\in  \X\times \X_g\}_{g\in G}.$    Notice that   $G^{(x,y)}=G^y,$ for all $x,y\in G$.  Then by (i) of Theorem \ref{VT} (applied to the partial action $\beta$) one has that
\[
x\in A_I \;\;\Longleftrightarrow \;\; (x,x)\in (\X^{2}\setminus A)^{*G}.
\]
Therefore $A_I$ is Borel thanks to Theorem \ref{VT3} (applied again to $\beta$).
\endproof

A theorem of Kechris  says that if $E$ is an equivalence relation on a Polish space, then $E$ has a Borel selector iff it is smooth and idealistic (\cite[Theorem 2.4]{KE1}, see also \cite[Theorem 5.4.11]{GA}). Using this fact  and Proposition \ref{idealisitic} we immediately obtain the following  generalization of a theorem of  Burgess   \cite[Corollary 5.4.12]{GA}.

\begin{teo}
\label{selector}
Let $m$ be a  partial action of $G$ on $\X.$ If $E^p_G$
is smooth, then it has a Borel selector.  
\end{teo}

\section{The enveloping space $\X_G$ is a standard Borel space}

In this section we prove that  the quotient Borel structure of $\X_G$ is standard, that is to say, we show that there is a topology $\tau$ on $\X_G$ such that $(\X_G, \tau)$ is Polish and  $Borel(\X_G, \tau)$ is the quotient Borel structure of $\X_G$. We also show that the orbit equivalence relation ${E}^p_G$ is Borel bireducible to an orbit equivalence relation induced by a total action of $G$ on some Polish space. 

The following is a generalization of (the easy part of)  a  theorem of  Effros \cite{EF}  (for more  details on this  topic see \cite{PU2}). 

\begin{teo}
\label{ef}
Let  $m$ be a topological partial action of  $G$ on $\X$.
 Then  the following assertions are equivalent.
\begin{enumerate}
\item[(i)]  $E_G^p$ is $G_\delta.$
\item[(ii)] $G^x\cdot x$ is $G_\delta$ in $\X,$ for every $x\in \X$.
\item[(iii)] $\X/E_G^p$ is $T_0.$ That is,  $ \overline{\{x\}}\neq \overline{\{y\}}$ for any $x,y\in\X/E_G^p $  with $x\ne y$.
\end{enumerate}
\end{teo}

\proof It is  shown exactly as in the case of a global action (see for instance \cite[Theorem 3.4.4]{GA}).
For (iii) implies (i), recall that  $\X/E_G^p$   is second countable by Lemma \ref{gxg} (ii). 
\endproof

\begin{lema}
\label{hat-selector}
Suppose that  $m$ is a   partial action of  $G$ on $\X.$
Then  $\widehat{E}^p_G$  is $G_\delta$ and has a Borel selector.
\end{lema}

\proof
If  $\widehat{E}^p_G$ is $G_\delta$, then it is  smooth thanks to Theorem \ref{gdelta} and thus   by Theorem \ref{selector} it has a Borel selector.
So it suffices to show that $ \widehat{E}^p_G$ is $G_\delta$ and by Theorem \ref{ef} we only need to show that  each orbit with respect to $\widehat m$ is $G_\delta.$  Let $(g,x)\in G\times \X$ and $O= G^{(g,x)}\cdot (g,x)$ the orbit of $(g,x)$.  Define $\rho: G^x\rightarrow  O$ by  $\rho(h)=(hg, h\cdot x)$.
Clearly  $\rho$ is a continuous bijection. Moreover,  $\rho^{-1}(j,y)=jg\m$,  then $\rho$ is a homeomorphism and therefore $O$ is Polish and hence $G_\delta$. 
\endproof

\begin{teo}
\label{gloto}
Let  $m$ be a  partial action of  $G$ on $\X$.   Then $\X_G$ is a $T_0$ space and its quotient Borel structure is standard. That is, there is  a Polish topology $\tau$ on $\X_G$ extending the given quotient topology of $\X_G$ such that  the quotient Borel structure of $\X_G$ is equal to $Borel(\X_G, \tau)$.  Moreover,   $Borel(\X)=Borel (\X,\tau_\X)$ and the enveloping action $\mu$ of $G$ on $\X_G$  is   $Borel(\X_G,\tau)$-measurable.
\end{teo}

\proof By Theorem \ref{mtomg}, $\X_G$ is equal to $(G\times \X)/\widehat{E}^p_G$. 
Thus by  Theorem \ref{ef} and Lemma \ref{hat-selector} one has that   $\X_G$ is $T_0.$ By Lemma \ref{hat-selector} there exists
$S:G\times\X\rightarrow G\times \X$   a Borel selector for $\widehat{E}^p_G$.  Then $T=\{ (g,x)\in G\times \X: \; S(g,x)=(g,x)\}$ is a transversal for  $\widehat{E}^p_G$, that is to say, $T$ has exactly one element in each $\widehat{E}^p_G$-equivalence class.   Consider the bijection $f: T\rightarrow (G\times\X)/ \widehat{E}^p_G$ given by $f(g,x)= \pi(g,x)$, where $\pi:G\times \X\rightarrow (G\times \X)/\widehat{E}^p_G$  is the quotient map.    We verify that $f$ is a Borel isomorphism.   It is clear that  $f$ is  continuous, thus Borel measurable.  Conversely, suppose $A\subseteq T$ is Borel, then by definition $\pi^{-1}(f[A])= S^{-1}(A)$ is also Borel.    Thus $f$ is a Borel isomorphism. Since by  \cite[Corollary 13.4]{KE} $T$  is a standard Borel space, then so is $(G\times\X)/ \widehat{E}^p_G$. 

We claim that w.l.o.g we can  assume that $S(1,x)=(1,x)$ for all $x\in\X$ and so $(1,x)\in T$ for all $x\in \X$. In fact, notice that  $(1,x)\widehat{E}^p_G (g,y)$ iff $g\in G^y$ and $g\cdot y=x$. Define $S'(g,x)= (1,g\cdot x)$ if $g\in G^x$ and $S'(g,x)=S(g,x)$ otherwise. Then $S'$ is also a Borel selector for $\widehat{E}^p_G$.

 Let $\tau$ be the topology on $\X_G$ induced by $f$, i.e., $V\subseteq \X_G$ is $\tau$-open iff $f\m (V)$ is open in $T$. Then $\tau$ extends the quotient topology on $\X_G$. It follows easily from the definition of $\iota$ (see \eqref{iota}) that $f\m(\iota (\X))=(G*\X)\cap T$   is a Borel subset of  $T$, then  $\X$ is $\tau$-Borel and $Borel(\X)=Borel (\X,\tau_\X)$.

Finally, since $\mu: G\times \X_G\rightarrow \X_G$ is continuous ($\X_G$ with the quotient topology), then it is clearly $Borel(\X_G,\tau)$-measurable.\endproof

The following example shows that the global action of $G$ on $\X_G$ is not necessarily continuous with respect to the topology given by Theorem \ref{gloto}.

\begin{exe}
\label{nonh}  \cite[Example 4.8]{PU}
Consider the partial action  of the discrete group  $\Z$ on $\mathbb{X}=[0,1]$ given by $m_{0}=\rm{id}_\mathbb{X},$ $m_{n}=\rm{id}_V,$ for $n\neq 0,$ and $V=(0,1].$  Notice that
$$
\Z*\mathbb{X}=\bigcup_{n\in \Z} \{n\}*\mathbb{X}=\{0\}* \mathbb{X}\,\,\cup\bigcup_{n\in \Z\atop n\neq 0} \{n\}*\mathbb{X}=\{0\}\times \mathbb{X} \,\cup\bigcup_{n\in \Z\atop n\neq 0} \{n\}\times V,
$$
is open in $\Z\times \mathbb{X}$. It is easy to verify that the equivalence classes of the form $[(n,0)]$ cannot be separated by open sets.  One may also notice that points in $\mathbb{X}_\Z$ are closed, so $\mathbb{X}_\Z$ is $T_1.$ In this case, by the proof of Theorem \ref{gloto}, the space $\X_\Z$ is Borel isomorphic to the Polish space $T=\{0\}\times [0,1]\,\cup \,\Z\times\{0\}.$ 
The enveloping action of $\Z$ on $\X_\Z$ defined in \eqref{action} is the following $\mu_m [(n,x)]=[(n+m,x)]$ for $n,m\in \Z$ and $x\in\X$. In particular, $\mu_1[(0,x)]=[(0,x)]$ for $x>0$, but $\mu_1[(0,0)] =[(1, 0)]$. 
Therefore $\mu$ is  not continuous at $[(0,0)]$.
\end{exe}

Since $\X_G$ is a standard Borel space, then it is natural to look at the orbit equivalence relation induced by the global action $\mu$. 
\begin{teo}
\label{bireducible}
Suppose that  $m$ is a   partial action of  $G$ on $\X.$  Then $E^p_G$ is Borel bireducible to the orbit equivalence relation induced by the enveloping  action of $G$ over $(\X_G,\tau)$, where $\tau$ is the topology given by Theorem \ref{gloto}.
\end{teo}

\proof Let $E_G$ denote the orbit equivalence relation induced $\mu$. It is immediate that $E^p_G\leq_B E_G$ as the map $\iota$ defined in \eqref{iota} is  a continuous embedding (with respect to the quotient topology on $\X_G$) and thus $Borel (\X_G,\tau)$-measurable.
Now we show that $E_G\leq_B E^p_G$.
By Lemma  \ref{hat-selector}, there exists $S:G\times \X\rightarrow G\times \X $ a Borel  selector for $\widehat{E}^p_G$. From the definition of $\widehat{E}_G^p$, we observe that  ${\rm proj}_2(S(g,x))E^p_Gx$, where ${\rm proj}_2:G\times \X\rightarrow \X$ is the projection to the second coordinate. Consider 
$$
f:\X_G\ni [g,x]\mapsto  {\rm proj}_2 ( S(g, x))\in \X.
$$ 
Since $S$ is a selector, then $f$ is well defined.  Let us see that $f$ is a Borel reduction between $E_G$ and $E^p_G$.  Clearly $f$ is Borel.
Suppose $[g,x]E_G[h,y]$, then there is $u\in G$ such that $ [ug, x]=[h,y]$. By Theorem \ref{mtomg},  one has that $(ug,x)\widehat{E}^p_G (h,y)$. Therefore, $S(ug,x)= S(h,y)$. Since
${\rm proj}_2 ( S(ug,x))  E_G^p   x$, then $ f([h,y])  E^p_G x$ and  $f([g,x])  E^p_G   x$. Thus $f([g,x])E^p_G f([h,y])$.
Conversely, if $ f([g,x])E^p_Gf([h,y]),$ then $xE^p_Gy,$ and there exists $\lambda \in G$ with $x\in \X_{\lambda\m}$ and $\lambda\cdot x=y.$ Let $u=h\lambda g\m,$ then $u[g,x]=[ug,x]=[h\lambda,x]=[h,\lambda\cdot x]=[h,y]$. Thus $[g,x]E_G[h,y],$ as desired. \endproof 

\end{document}